\documentclass{elsarticle}
\journal{Linear Algebra and its Applications}
\usepackage[utf8]{inputenc}
\usepackage[T1]{fontenc}
\usepackage{amsfonts,amsmath,amssymb,amsthm}

\usepackage{hyperref}
\usepackage{pgfplots}
\usepackage{algorithmic,algorithm}
\usepackage[colorinlistoftodos,textwidth=4cm,shadow]{todonotes}
\usepackage{multirow}
\newtheorem{theorem}{Theorem}[section]
\newdefinition{definition}{Definition}
\newtheorem{remark}[theorem]{Remark}
\newtheorem{lemma}[theorem]{Lemma}
\newcounter{Ivan}

\newcounter{Sasha}

\begin{document}


\begin{frontmatter}

\title{Rectangular maximum-volume submatrices and their applications}
\author[KAUST]{A.~Mikhalev\corref{corauthor}}
\ead{aleksandr.mikhalev@kaust.edu.sa}
\ead{muxasizhevsk@gmail.com}
\author[SKOLTECH,INMRAS]{I.~V.~Oseledets}
\ead{i.oseledets@skoltech.ru}
\address[KAUST]{King Abdullah University of Science and Technology,
Thuwal 23955-6900, Kingdom of Saudi Arabia.}
\address[SKOLTECH]{Skolkovo Institute of Science and Technology,
Novaya St.~100, Skolkovo, Odintsovsky district, 143025, Russia.}
\address[INMRAS]{Institute of Numerical Mathematics, Russian Academy of
Sciences. Gubkina~St.~8, 119333 Moscow, Russia.}
\cortext[corauthor]{Corresponding author}


\begin{abstract}
We introduce a definition of the volume of a general rectangular matrix,
which is equivalent to an absolute value of the determinant for square
matrices. We generalize results of square maximum-volume submatrices to the
rectangular case, show a connection of the rectangular volume with an optimal
experimental design and provide estimates for a growth of coefficients and an
approximation error in spectral and Chebyshev norms. Three promising
applications of such submatrices are presented: recommender systems, finding
maximal elements in low-rank matrices and preconditioning of overdetermined
linear systems. The code is available online.
\end{abstract}



\end{frontmatter}


\section{Introduction}
\label{sec:Intro}
How to define the volume of a rectangular matrix, and how to compute a
submatrix with the maximal volume in a given matrix? A standard definition of
the volume of a square matrix is an absolute value of its determinant.
Maximum-volume submatrices play an important role in low-rank
approximations \cite{gtz-psa-1997,gt-maxvol-2001}, recommender systems
\cite{liu2011wisdom}, wireless communications \cite{wang2010global},
preconditioning of overdetermined systems \cite{arioli2015preconditioning},
tensor decompositions \cite{ot-ttcross-2010}. 
How to compute a submatrix of exactly maximal volume is a NP-hard problem
\cite{bartholdi-1982}. However, in many applications, a submatrix of a
sufficiently large volume is enough, and it can be computed in a polynomial
time using the \textbf{maxvol} algorithm \cite{gostz-maxvol-2010}. The
\textbf{maxvol} algorithm is a greedy iterative algorithm, which swaps rows to
maximize the volume of a square submatrix.

In this paper, we extend the volume concept to the case of \emph{rectangular
matrices} (section \ref{sec:Volume}), introduce a \emph{dominance} property,
which is important for theoretical estimations (section \ref{sec:dominant}),
generalize well-known results of the square case (sections \ref{sec:bounds}
and \ref{sec:norm}), remind of pseudo-skeleton and CGR approximations and
provide estimations of approximation error (section \ref{sec:approx}),
propose new volume maximization algorithm (so-called \textbf{rect\_maxvol},
section \ref{sec:rect_maxvol}) and apply the
\textbf{rect\_maxvol} algorithm for three different problems and compare its
results with results of the \textbf{maxvol} algorithm (section
\ref{sec:numerical}). We also show a connection of a new definition of the
volume with an optimal experimental design (section \ref{sec:optdesign}).

\section{Volume of rectangular matrices}
\label{sec:Volume}

The volume of a square matrix $A$ has a natural geometric meaning as a volume
of the parallelepiped, spanned by rows of the matrix $A$,
and is equal to the product of its singular values. This
definition can be straightforwardly generalized to the rectangular case as
$\sqrt{\det A^*A}$ or $\sqrt{\det AA^*}$, depending on the shape of $A$.
Let us assume, that a number of rows of the matrix $A$ is not less than a
number of columns. So, we use $\sqrt{\det A^*A}$ as a value of the rectangular
volume of the matrix $A$:
\begin{equation*}
\mathrm{vol}(A) = \sqrt{\det(A^*A)}.
\end{equation*}

Geometric meaning of this definition is the following: for a $K$-by-$r$ matrix
$A$
with $K \ge r$, it shows how many times an $r$-dimensional Euclidean volume of
an image of a unit $K$-dimensional ball under a linear operator $A^*$ is
greater than Euclidean volume of a unit $r$-dimensional ball. This fact has a
rather simple proof.
Following a singular values decomposition, matrix $A^*$ can be combined as
a multiplication of an $r$-by-$r$ unitary, an $r$-by-$K$ diagonal and a
$K$-by-$K$ unitary matrices.
A linear operator $A^*$ does the following changes to a unit ball: rotate
(which
does not change a unit ball at all), project to an $r$-dimensional space
(leaving us with an $r$-dimensional unit ball), scale axes by singular values
and rotate in an $r$-dimensional space. Since a rotation does not change
Euclidean volume, the only step, that changes it, is scaling axes by singular
values. As the product of singular values of the matrix $A^*$ is equal to
$\sqrt{\det(A^*A)}$, we finish the proof.

\section{Dominance property}
\label{sec:dominant}
As it was already mentioned, it is a NP-hard problem to find the exact
maximal-volume submatrix. That is why in order to find a good submatrix
in a reasonable amount of time, the maximal-volume property is typically
relaxed to a so-called \emph{dominance}
\cite{gtz-psa-1997,gt-maxvol-2001,gostz-maxvol-2010} property.
Standard definition of a square \emph{dominant} submatrix is the following:
\begin{definition}[Square dominant submatrix]
\label{def:dominant}
Let $N \ge r$ and $A \in \mathbb{C}^{N \times r}$ be of a full column rank. An
$r \times r$ submatrix $\widehat{A}$ is called a \emph{dominant} submatrix, if
a swap of any single row of $\widehat{A}$ for a row of $A$, not already
presented in $\widehat{A}$, does not increase the volume.
\end{definition}
For a simplicity, we assume that $\widehat{A}$ corresponds to an upper part of
$A$ and a complementary submatrix $\widetilde{A}$ corresponds to a lower part
of $A$:
\begin{equation*}
A = \begin{bmatrix} \widehat{A} \\ \widetilde{A} \end{bmatrix}.
\end{equation*}
Let us denote $C$ as the following matrix of coefficients:
\begin{equation*}
C = A \widehat{A}^{-1} = \begin{bmatrix} I_{r \times r} \\ \widetilde{C}
\end{bmatrix}.
\end{equation*}
In \cite{gostz-maxvol-2010}, it was shown, that a swap of an $i$-th row of
$\widehat{A}$ by a $j$-th row of $\widetilde{A}$ multiplies the volume of
$\widehat{A}$ by a modulus of $\widetilde{C}_{ji}$, an element of
$\widetilde{C}$ on intersection of a $j$-th row and an $i$-th column. If
$\widehat{A}$ is a dominant submatrix, then elements of the matrix
$\widetilde{C}$ are less than or equal to $1$ in modulus. So, the Chebyshev
($l_\infty$) norm of the matrix $C$ is bounded:
\begin{equation}
\label{dominant:C}
\Vert C \Vert_C \le 1.
\end{equation}
A geometric meaning of \eqref{dominant:C} is that among any $N$ vectors in an
$r$-dimensional space we can select $r$ vectors in a such way, that any vector
out of given $N$ vectors can be expressed as a linear combination of the
selected vectors with coefficients, less than $1$ in modulus. The inequality
\eqref{dominant:C} can not be improved if we consider only square submatrices.
In a practice, it is much easier to find a \emph{quasi-dominant} submatrix,
which is only approximately \emph{dominant}:
\begin{definition}[Square quasi-dominant submatrix]
Let $N \ge r$ and $A \in \mathbb{C}^{N \times r}$ be of a full column rank. An
$r \times r$ submatrix $\widehat{A}$ is called a \emph{quasi-dominant}
submatrix, if a
swap of any single row of $\widehat{A}$ for a row of $A$, not already presented
in $\widehat{A}$, does not increase the volume by more than a factor of
$1+\varepsilon$.
\end{definition}
The \textbf{maxvol} \cite{gostz-maxvol-2010} algorithm, a baseline method for
our numerical experiments, finds a \emph{quasi-dominant} submatrix. The matrix
of coefficients, corresponding to a \emph{quasi-dominant} submatrix, is bounded
as
\begin{equation*}
\Vert C \Vert_C \le 1+\varepsilon
\end{equation*}

Definitions of \emph{dominant} and \emph{quasi-dominant} submatrices can be
naturally extended to the rectangular case:
\begin{definition}[Rectangular dominant submatrix]
Let $N \ge K \ge r$ and $A \in \mathbb{C}^{N \times r}$ be of a full column
rank. A rectangular $K \times r$ submatrix $\widehat{A}$ is called a
\emph{dominant} submatrix, if a swap of any single row of $\widehat{A}$ by a
row of $A$,
not already presented in $\widehat{A}$, does not increase the volume.
\end{definition}
\begin{definition}[Rectangular quasi-dominant submatrix]
Let $N \ge K \ge r$ and $A \in \mathbb{C}^{N \times r}$ be of a full column
rank. A rectangular $K \times r$ submatrix $\widehat{A}$ is called a
\emph{quasi-dominant} submatrix, if a swap of any single row of $\widehat{A}$
by a row of
$A$, not already presented in $\widehat{A}$, does not increase the volume by
more, than a factor of $1+\varepsilon$.
\end{definition}

The \emph{dominance} property plays important role in theoretical estimations,
presented later in this paper (section~\ref{sec:bounds}). However, submatrices,
selected by our \textbf{rect\_maxvol} method (presented in section
\ref{sec:rect_maxvol}) are not even \emph{quasi-dominant}
by a construction. Nevertheless, numerical experiments
(section~\ref{sec:numerical}) show promising results.

\section{Main result}

In this section, we derive theoretical properties of rectangular dominant and
maximal-volume submatrices (section \ref{sec:bounds}), provide a spectral
analysis of arising matrices of coefficients (section \ref{sec:norm}), propose
special constructions for pseudo-skeleton and CGR approximations and show
their influence on approximation errors in Chebyshev and spectral norms
(section \ref{sec:approx}) and propose the \textbf{rext\_maxvol} algorithm for
finding ``good'' rectangular submatrices (section \ref{sec:rect_maxvol}).

Note that we follow the philosophy of the paper \cite{gtz-psa-1997} for the
square case.

\subsection{Upper bound on coefficients}
\label{sec:bounds}
We start with a simple lemma, followed by a theorem on upper bounds for a
matrix of coefficients:
\begin{lemma}
\label{rectmaxvol:detlem}
Let $M > N$, a matrix $A \in \mathbb{C}^{N \times M}$ and a matrix $B \in
\mathbb{C}^{M \times N}$. Let $A_{-i}$ be a $N \times (M-1)$ submatrix of $A$
without $i$-th column and $B_{-i}$ be a $(M-1) \times N$ submatrix of $B$
without $i$-th row. Then,
\begin{equation*}
\det(AB) = \frac{1}{M-N}\sum_{i=1}^M\det(A_{-i}B_{-i}).
\end{equation*}
\end{lemma}
\begin{proof}
From the Cauchy-Binet formula we get:
\begin{equation*}
\det(AB) = \sum_{(j)} \det(A_{(j)})\det(B_{(j)}),
\end{equation*}
where $(j)$ is a set of $N$ different numbers, such that $A_{(j)}$ is a
submatrix on columns $(j)$ of the matrix $A$ and $B_{(j)}$ is a submatrix on
rows $(j)$ of the matrix $B$. Let us assume any single set $(j)$. A
submatrix $A_{-i}$ contains all columns of $A$, except $i$-th column, so
$A_{(j)}$ is a submatrix of $A_{-i}$ for any $i \notin (j)$. Since $(j)$
consists of $N$ different numbers and $M$ is a total number of columns of $A$,
there are $M-N$ different values of $i$ such that $A_{(j)}$ is a submatrix of
$A_{-i}$. The same goes for the submatrix $B_{(j)}$. So, according to the
Cauchy-Binet formula, multiplication of determinants $\det{A_{(j)}}
\det{B_{(j)}}$ is a summand of $\det(A_{-i}B_{-i})$ if and only if $i \notin
(j)$ with $M-N$ possible values of $i$. So, we get
\begin{equation*}
\sum_{i=1}^M \det(A_{-i}B_{-i}) = (M-N)\sum_{(j)} \det(A_{(j)})\det(B_{(j)}) =
(M-N)\det(AB),
\end{equation*}
and finish the proof.
\end{proof}
\begin{theorem}
\label{rectmaxvol:est}
Let $N \ge K \ge r$ and a matrix $A \in \mathbb{C}^{N \times r}$ be of a rank
$r$. Let $\widehat{A}$ be a $K \times r$ dominant submatrix of the matrix $A$,
based on a set of rows $(j)$. Then, there is such matrix of coefficients $C$,
that $A = C \widehat{A}$ and for every row, excluding rows from the set $(j)$,
the following inequality holds:
\begin{equation*}
\forall i \in \{1, \ldots, N\} \setminus (j):\:\Vert C_i \Vert_2 \le \sqrt{
\frac{r}{K+1-r}}.
\end{equation*}
\end{theorem}
\begin{proof}
Since the matrix $A$ is of a full column rank, its dominant $K$-by-$r$
submatrix is non-singular. It means, that for every $i$ solution $C_i$ of an
equation
\begin{equation*}
C_i \widehat{A} = A_i
\end{equation*}
exists, but it may be not unique. Let us use just any solution in the case $i
\in (j)$ and the minimum norm solution in the case $i \notin (j)$. So, $C$ is a
solution of
\begin{equation*}
C \widehat{A} = A.
\end{equation*}
Let us assume $i \notin (j)$ and construct a matrix $H$ as follows:
\begin{equation*}
H = \begin{bmatrix} \widehat{A} \\ A_i \end{bmatrix}.
\end{equation*}
From the determinant equation for the Sch\"{u}r complement:
\begin{equation*}
\det(H^* H) = \det(\widehat{A}^* \widehat{A} + A_i^* A_i) = \det(\widehat{A}^*
\widehat{A}) (1 + A_i (\widehat{A}^* \widehat{A})^{-1} A_i^*).
\end{equation*}
Since $C_i$ is the minimum norm solution,
\begin{equation*}
A_i (\widehat{A}^* \widehat{A})^{-1} A_i^* = C_i C_i^*,
\end{equation*}
\begin{equation}
\label{eq:volume_grow}
\det(H^* H) = \det(\widehat{A}^* \widehat{A}) (1+\Vert C_i \Vert_2^2),
\end{equation}
\begin{equation*}
\Vert C_i \Vert_2^2 = \frac{\det(H^* H)}{\det(\widehat{A}^* \widehat{A})}-1.
\end{equation*}
The submatrix $\widehat{A}$ is a dominant $K \times r$ submatrix of the matrix
$A$, so it has maximum volume among all $K \times r$ submatrices of the
matrix $H$. Applying lemma~\ref{rectmaxvol:detlem} to the matrix $H$, we get
\begin{equation*}
\det(H^* H) \le \frac{K+1}{K+1-r}\det(\widehat{A}^* \widehat{A}).
\end{equation*}
So, we have an upper bound on $l_2$ norm of the $i$-th row of $C$:
\begin{equation*}
\Vert C_i \Vert^2_2 \le \frac{r}{K+1-r}.
\end{equation*}
Applying the latter inequality for every $i$ not in the set $(j)$ we complete
the proof.
\end{proof}

Note that the same result of theorem~\ref{rectmaxvol:est} was obtained in
\cite{de2007subset}, where it was used to estimate the Frobenius norm of the
matrix of coefficients $C$. However, we put here our own proof since it
provides an interesting alternative. Our main goal here is to show that a
dominant submatrix leads to an upper bound for rows of the matrix of
coefficients. Later, in section~\ref{sec:approx}, we will show, that the
Chebyshev norm of an approximation error practically linearly depends on this
bound.

\subsection{Spectral analysis of two different matrices of coefficients}
\label{sec:norm}
Without the loss of generality, let the dominant submatrix $\widehat{A} \in
\mathbb{C}^{K \times r}$ be located in the first rows of $A \in \mathbb{C}^{N
\times r}$:
\begin{equation*}
A = \begin{bmatrix} \widehat{A} \\ B \end{bmatrix}.
\end{equation*}
So, the matrix of coefficients $C \in \mathbb{C}^{N \times K}$ can be divided
into submatrices:
\begin{equation*}
C = \begin{bmatrix} \widehat{C} \\ B \widehat{A}^\dagger \end{bmatrix}.
\end{equation*}
From the theorem~\ref{rectmaxvol:est} we got the bound for rows of the
matrix $B \widehat{A}^\dagger$, by putting those rows equal to the
minimum norm solutions of corresponding equations with the matrix
$\widehat{A}$. However, in the rectangular case $\widehat{C} \in \mathbb{C}^{K
\times K}$ is not unique and we have to set it to a some reasonable value. Two
obvious variants are $\widehat{C} = I_{K \times K}$ and
$\widehat{C} = \widehat{A} \widehat{A}^\dagger$. In this section we show
singular values of $C$ for both variants.

We need to show that singular values of following matrices are
practically the same:
\begin{equation*}
C_1 = \begin{bmatrix} I_{K \times K} \\ B \widehat{A}^\dagger \end{bmatrix},
\;\; C_2 = \begin{bmatrix} \widehat{A} \widehat{A}^\dagger \\ B
\widehat{A}^\dagger \end{bmatrix} = C_1 \widehat{A} \widehat{A}^\dagger.
\end{equation*}
First of all, singular values of $C_1$ are, obviously, following:
\begin{equation*}
\sigma_i(C_1) = \sqrt{1+\sigma_i^2(B \widehat{A}^\dagger)}.
\end{equation*}
Since $\widehat{A} \widehat{A}^\dagger$ is an orthoprojector to the space of
the first $r$ right singular vectors of $\widehat{A}^\dagger$, the first $r$
singular values of $C_2$ are equal to the first $r$ singular values of $C_1$,
while all other singular values of $C_2$ are zero. So, we get following
equations for singular values:
\begin{equation}
\label{eq:topsingval}
\forall i=1..r: \sigma_i(C_1) = \sigma_i(C_2) = \sqrt{1+\sigma_i^2(B
\widehat{A}^\dagger)},
\end{equation}
\begin{equation*}
\forall i=(r+1)..K: \sigma_i(C_1) = 1, \sigma_i(C_2) = 0.
\end{equation*}

As it can be seen, the spectral norm of the matrix $C$ does not depend on
a selection of the submatrix $\widehat{C}$ in examined cases. However,
$\widehat{C} = I_{K \times K}$ is more intuitive to use.

\subsection{Rectangular pseudo-skeleton and CGR-approximations}
\label{sec:approx}
Skeleton type approximations of a given matrix $A$ are based on specially
selected rows $R$, columns $C$ and a core matrix $G$:
\begin{equation}
\label{eq:CGR}
A \approx C G R.
\end{equation}
For a CGR-approximation (also known as CUR-approximation), a core matrix $G$
can
be chosen in any convenient way, while for a pseudo-skeleton approximation, a
core matrix is a pseudo-inverse of a submatrix on an intersection of the
rows $R$ and the columns $C$. Error estimations in the case of equal number of
specially selected rows and columns can be found in \cite{gtz-psa-1997} and
\cite{gt-maxvol-2001}. To estimate the error in the case of rectangular
pseudo-skeleton or CGR approximation in the spectral norm, we need to remind a
definition from \cite{gtz-psa-1997} and add an additional one:

\begin{definition}[$t(r,n)$ \cite{gtz-psa-1997}]
\label{def:trn}
Let $n \ge r$. Let $\mathcal{P}(n, r)$ be a space of all $n \times r$
orthogonal matrices (Stiefel manifold \cite{james1976topology}). Let denote
$\mathcal{M}(U)$ as a set of all $r \times r$ submatrices of a given
orthogonal matrix $U$ and $\sigma_{min}(\widehat{U})$ as the minimal singular
value of a matrix $\widehat{U}$. Then, define $t(r, n)$ as follows:
\begin{equation*}
t(r, n) = \left[ \min_{U \in \mathcal{P}(n,r)} \left( \max_{\widehat{U} \in
\mathcal{M}(U)} \sigma_{\min} (\widehat{U}) \right) \right]^{-1}.
\end{equation*}
\end{definition}
\begin{definition}[$t(r,n,k)$]
\label{def:trnk}
Let $n \ge k \ge r$. Let $\mathcal{P}(n,r )$ be a space of all $n \times r$
orthogonal matrices (Stiefel manifold \cite{james1976topology}). Let denote
$\mathcal{M}_k(U)$ as a set of all $k \times r$ submatrices of a given
orthogonal matrix $U$ and $\sigma_{min}(\widehat{U})$ as the minimal singular
value of a matrix $\widehat{U}$. Then, define $t(r, n, k)$ as follows:
\begin{equation*}
t(r, n ,k) = \left[ \min_{U \in \mathcal{P}(n,r)} \left( \max_{\widehat{U}
\in \mathcal{M}_k(U)} \sigma_{\min} (\widehat{U}) \right) \right]^{-1}.
\end{equation*}
\end{definition}

One can show, that the inner maximum of definitions \ref{def:trn} and
\ref{def:trnk} over all submatrices is a continuous
function of $U$. Since Stiefel manifold is compact and the inner maximum is a
continuous function, the outer minimum is achievable on a some orthogonal
matrix.
So, there is a such orthogonal matrix $U$ with a such submatrix $\widehat{U}$,
that
\begin{equation*}
t(r, n, k) = \frac{1}{\sigma_{min}(\widehat{U})}.
\end{equation*}
Definition \ref{def:trnk} is a formal generalization of $t(r, n)$, described in
\cite{gtz-psa-1997}, to the case of rectangular submatrices. The meaning of the
$t(r, n, k)$ is very simple: any $n \times r$ orthogonal matrix has such $k
\times r$ submatrix, that norm of the pseudo-inverse of this submatrix is
upper-bounded by $t(r, n, k)$.

\begin{lemma}
\label{lemma:trnk}
For any given $r \le k \le n$, value $t(r, n, k)$ has the following upper
bound:
\begin{equation*}
t(r, n, k) \le \sqrt{1+\frac{(n-k)r}{k+1-r}}.
\end{equation*}
\end{lemma}

\begin{proof}
In the definition \ref{def:trnk} instead of the inner maximum over all
submatrices of the matrix $U$ we can use any dominant
submatrix. Let $U$ be orthogonal and $\widehat{U}$ be its
dominant submatrix. Let $\widetilde{U}$ be a submatrix, complementary to
dominant. From equation \eqref{eq:topsingval} we get the spectral norm of
$\widehat{U}^\dagger$:
\begin{equation*}
t(r, n, k) \le \Vert \widehat{U}^\dagger \Vert_2 = \Vert U \widehat{U}^\dagger
\Vert_2 = \Vert C \Vert_2 = \sqrt{1 + \Vert \widetilde{U} \widehat{U}^\dagger
\Vert_2^2}.
\end{equation*}
The spectral norm of $\widetilde{C} = \widetilde{U} \widehat{U}^\dagger$ is
upper
bounded by its Frobenius norm, which can be bounded by the
theorem~\ref{rectmaxvol:est}:
\begin{equation*}
\Vert \widetilde{C} \Vert_F^2 = \frac{(n-k)r}{k+1-r}.
\end{equation*}
So, we got upper estimates for $t(r, n, k)$.
\end{proof}
Now we can check several values of $k$:
\begin{equation*}
t(r, n) = t(r, n, r) \le \sqrt{(n-r)r+1},
\end{equation*}
\begin{equation*}
t(r, n, 1.25r-1) \le \sqrt{4n-5r+5},
\end{equation*}
\begin{equation*}
t(r, n, 2r-1) \le \sqrt{n-2r+2}.
\end{equation*}

Authors of \cite{gtz-psa-1997} proposed the hypothesis \footnote{This
hypothesis is not yet proven}
\begin{equation*}
t(r, n) \le \sqrt{n},
\end{equation*}
%
Unfortunately, we were not able to provide a similar hypothesis for
$t(r, n, k)$.

Now we proceed to error estimations of rectangular pseudo-skeleton
approximations. We use specially constructed approximants to prove upper bounds
of an approximation error. These approximants use so-called ``basis'' rows and
columns, which have a very simple definition:

\begin{definition}[``basis'' rows (columns)]
\label{def:basisrows}
Let $A$ be a $N$-by-$r$ ($r$-by-$N$) matrix of a full column (row) rank. Then,
given $n \ge r$ rows (columns) are called ``basis'' if any other row (column)
of $A$ can be written as a linear combination of given ones.
\end{definition}

Of course, one has to be careful when selecting such ``basis'' rows or
columns, since it influences the overall approximation error directly, which can
be
seen in proofs of theorems \ref{thm:pseudoapproximation},
\ref{thm:CGRapproximation} and \ref{thm:CGRapproximation2}.
As we provide estimations in spectral (theorems \ref{thm:pseudoapproximation}
and \ref{thm:CGRapproximation}) and Chebyshev (theorem
\ref{thm:CGRapproximation2}) norms, we propose to select ``basis'' rows differently as in the following remark.
\begin{remark}[How to select ``basis'' rows]
\label{remark:basis}
Let a matrix $A$ be a $N$-by-$r$ matrix of a full column rank, a matrix
$\widehat{A}$
be a submatrix on ``basis'' rows and matrix $C$ be the minimal norm solution of
\begin{equation*}
C \widehat{A} = A.
\end{equation*}
Let us select ``basis'' rows (``basis'' submatrix) of $A$ as
\begin{enumerate}
\item a dominant $n$-by-$r$ submatrix of $A$ in the case of
estimations in the Chebyshev norm, with
$$\Vert C_i \Vert_2 \le \sqrt{\frac{r}{n+1-r}}$$ for each non-``basis'' row
$C_i$ of the matrix $C$,
\item a $n$-by-$r$ submatrix with the minimal possible
spectral norm of the corresponding matrix $C$ in the case of estimations in
the spectral norm, with
$$\Vert C \Vert_2 \le t(r, N, n).$$
\end{enumerate}
\end{remark}

The selection technique from this remark has a one drawback: it is a NP-hard
problem
to acquire such
subsets in both spectral and Chebyshev cases. A practical way to select such
``basis'' is given only for the case of the Chebyshev norm and is presented in
the next
section \ref{sec:rect_maxvol}, but it is convenient to use the remark
\ref{remark:basis} for theoretical proofs.

So, we define our approximants:

\begin{definition}[Rectangular pseudo-skeleton approximant]
\label{def:pseudoskeleton}
Let $A \in \mathbb{C}^{N \times M}$, $Z = Z_U Z_V$ be its low-rank
approximation with $Z_U \in \mathbb{C}^{N \times r}$ and $Z_V \in \mathbb{C}^{r
\times M}$. Then, construct a pseudo-skeleton approximant by following steps:
\begin{enumerate}
\item Define $n$ ``basis'' rows of $Z_U$,
\item Denote corresponding rows of $A$ as a matrix $R$,
\item Define $m \ge n$ ``basis'' columns of $R$,
\item Denote corresponding columns of $A$ as a matrix $C$,
\item Denote $\widehat{A}$ as a submatrix on intersection of rows $R$ and
columns $C$,
\item $C \widehat{A}^\dagger R$ is a rectangular pseudo-skeleton approximant.
\end{enumerate}
\end{definition}

\begin{definition}[Rectangular CGR-approximant]
\label{def:CGR}
Let a matrix $A$ be $N$-by-$M$ complex or real matrix, a matrix $Z$ be its best
rank-$r$ approximation in spectral and Frobenius norms. Then, we propose
following steps to construct a CGR-approximant:
\begin{enumerate}
\item Factorize $Z = Z_U Z_V$ with $Z_U \in \mathbb{C}^{N \times r}$ and $Z_V
\in \mathbb{C}^{r \times M}$,
\item Define $n$ ``basis'' rows of $Z_U$,
\item Denote corresponding rows of $A$ as a matrix $R$,
\item Compute a singular values decomposition of $R$, truncate singular values
and vectors, starting from $(r+1)$-th, and denote $USV$ as a main part and $E$
as a truncated part or a noise,
\item Define $m$ ``basis'' columns of $V$,
\item Denote corresponding submatrix of $V$ as $W$,
\item Denote corresponding columns of $A$ as a matrix $C$,
\item Define kernel matrix $G$ as a matrix $(USW)^\dagger$,
\item $C (USW)^\dagger R$ is a rectangular CGR-approximant.
\end{enumerate}
\end{definition}

Proposed definitions are correspondingly used in the following theorems.

\begin{theorem}[Rectangular pseudo-skeleton approximation error]
\label{thm:pseudoapproximation}
Let a matrix $A \in \mathbb{C}^{N \times M}, A=Z+F, \mathrm{rank}\,Z = r,
\Vert F \Vert \le \varepsilon$. Then, the error of approximation by the
rectangular
pseudo-skeleton approximant $\widetilde{A}$ (from the definition
\ref{def:pseudoskeleton}), based on $n$ rows and $m$ columns $(n \ll N, m \ll
M, m \ge n \ge r)$, has the following upper bound:
\begin{equation*}
\Vert A - \widetilde{A} \Vert_2 \le \varepsilon t(n, M, m) t(r, N, n).
\end{equation*}
\end{theorem}
\begin{proof}
Assume matrices $A, Z$ and $F$ are divided into blocks
\begin{equation*}
A = \begin{bmatrix}A_{11}&A_{12}\\A_{21}&A_{22}\end{bmatrix}, \;\;
Z = \begin{bmatrix}Z_{11}&Z_{12}\\Z_{21}&Z_{22}\end{bmatrix}, \;\;
F = \begin{bmatrix}F_{11}&F_{12}\\F_{21}&F_{22}\end{bmatrix}.
\end{equation*}
Without the loss of generality, let the submatrix $A_{11} \in \mathbb{C}^{n
\times m}$ be a core matrix for the pseudo-skeleton approximation
\eqref{eq:CGR}:
\begin{equation}
\label{eq:pseudoskel}
\widetilde{A} = \begin{bmatrix}A_{11}\\A_{21}\end{bmatrix} A_{11}^\dagger
\begin{bmatrix} A_{11} & A_{12}\end{bmatrix}.
\end{equation}
Since $A_{11}$ was chosen as a ``basis'' submatrix of $\begin{bmatrix}A_{11}&
A_{12}\end{bmatrix}$, a matrix $A_{11}^\dagger \begin{bmatrix}
A_{11}&A_{12}\end{bmatrix}$ is well-defined, and the total pseudo-skeleton 
approximation \eqref{eq:pseudoskel} makes sense.

Let a matrix $C_Z$ be the minimum norm solution of
\begin{equation}
\label{eq:C_Z}
\begin{bmatrix}Z_{21}&Z_{22}\end{bmatrix} = C_Z
\begin{bmatrix}Z_{11}&Z_{12}\end{bmatrix},
\end{equation}
and $C_A$ be the minimum norm solution of
\begin{equation*}
A_{12} = A_{11} C_A.
\end{equation*}
By construction,
\begin{equation}
\label{eq:approx0}
\widetilde{A}-A = \begin{bmatrix} 0_{n \times m} & 0_{n \times (M-m)}
\\A_{21}A_{11}^\dagger A_{11}-A_{21}&
A_{21}A_{11}^\dagger A_{12}-A_{22}\end{bmatrix}.
\end{equation}
Let us estimate the value $A_{21}A_{11}^\dagger A_{11}$:
\begin{multline*}
A_{21}A_{11}^\dagger A_{11} = (Z_{21}+F_{21})A_{11}^\dagger A_{11} =
(C_Z (A_{11}-F_{11})+F_{21})A_{11}^\dagger A_{11} = \\ =
C_Z(Z_{11}+F_{11})+(F_{21}-C_Z F_{11})A_{11}^\dagger A_{11} = \\ =
A_{21}-(F_{21}-C_Z F_{11})(I_{m \times m}-A_{11}^\dagger A_{11}),
\end{multline*}
and its difference with $A_{21}$:
\begin{equation}
\label{eq:A21}
A_{21}A_{11}^\dagger A_{11}-A_{21} = \begin{bmatrix}C_Z & -I_{(N-n) \times
(N-n)} \end{bmatrix} F \begin{bmatrix}I_{m \times m}-A_{11}^\dagger A_{11}
\\ 0_{(M-m) \times m} \end{bmatrix}.
\end{equation}
Let us also estimate an approximation of $A_{22}$:
\begin{multline*}
A_{21}A_{11}^\dagger A_{12} = (Z_{21}+F_{21})C_A = (C_Z
(A_{11}-F_{11})+F_{21})C_A = \\ =
C_Z A_{11} C_A+(F_{21}-C_Z F_{11})C_A,
\end{multline*}
submatrix $A_{22}$ itself:
\begin{equation*}
A_{22} = C_Z Z_{12}+F_{22} = C_Z
(A_{11} C_A-F_{12})+F_{22} = C_Z A_{11}
C_A+(F_{22}-C_Z F_{21}),
\end{equation*}
and their difference:
\begin{multline}
\label{eq:A22}
A_{21}A_{11}^\dagger A_{12}-A_{22} = (F_{21}-C_Z F_{11})C_A-
(F_{22}-C_Z F_{21}) = \\ =\begin{bmatrix}C_Z & -I_{(N-n) \times (N-n)}
\end{bmatrix} F \begin{bmatrix}-C_A \\ I_{(M-m) \times (M-m)} \end{bmatrix}.
\end{multline}
Combining equations \eqref{eq:approx0}, \eqref{eq:A21} and \eqref{eq:A22},
we get
\begin{equation}
\label{eq:approx}
\widetilde{A}-A = L F R,
\end{equation}
where
\begin{equation*}
L = \begin{bmatrix} 0_{n \times n} & 0_{n \times (N-n)} \\ C_Z &
-I_{(N-n) \times (N-n)}\end{bmatrix}, \;\;
R = \begin{bmatrix} I_{m \times m}-A_{11}^\dagger A_{11} & -C_A \\ 0_{(M-m)
\times m} & I_{(M-m) \times (M-m)}\end{bmatrix}.
\end{equation*}
Obviously,
\begin{equation}
\label{eq:L}
\left \Vert L \right \Vert_2 = \sqrt{1+\Vert C_Z\Vert_2^2}.
\end{equation}
The matrix $I_{m \times m}-A_{11}^\dagger A_{11}$ is symmetric and orthogonal
to matrix $C_A$, so the first $m$ columns of $R$ are orthogonal to all other
columns of $R$. Since the spectral norm of the first $m$ columns of $R$ is $1$
or $0$, depending on a relation of $n$ to $m$, and the spectral norm of other
columns is not less than 1, we get
\begin{equation}
\label{eq:R}
\left \Vert R \right \Vert_2 = \sqrt{1+\Vert C_A \Vert_2^2}.
\end{equation}
Matrix $Z$ is of rank $r$ and $C_Z$ is the minimal norm solution of
\eqref{eq:C_Z}, so, due to Remark \ref{remark:basis},
%
we have the following upper bound:
\begin{equation}
\label{eq:CZ}
\sqrt{1+\Vert C_Z \Vert_2^2} \le t(r, N, n).
\end{equation}
Using the same technique for the $Q$ factor of the QR-factorization of
$\begin{bmatrix} A_{11} & A_{12} \end{bmatrix}$, we get
\begin{equation}
\label{eq:CA}
\sqrt{1+\Vert C_A \Vert_2^2} \le t(n, M, m).
\end{equation}
Combining equations \eqref{eq:approx}, \eqref{eq:L}, \eqref{eq:R},
\eqref{eq:CZ} and \eqref{eq:CA} we finish the proof for an error estimation in
the
spectral norm.
\end{proof}

Estimations, provided in theorem~\ref{thm:pseudoapproximation}, are not
symmetric due to the construction of an approximation, which can be slightly
changed
to make estimations symmetric. However, this changes a pseudo-skeleton
approximation to a CGR-approximation \eqref{eq:CGR} with a specially selected
kernel matrix.

\begin{theorem}[Rectangular CGR-approximation error]
\label{thm:CGRapproximation}
Let  matrix $A$ be an $N$-by-$M$ complex or real matrix, matrix $Z$ be its best
rank-$r$ approximation in spectral and Frobenius norms. Then, a rectangular
CGR-approximant $\widetilde{A}$ (from Definition \ref{def:CGR}), based on
$n$
rows ($r \le n \ll N$) and $m$ columns ($r \le m \ll M$), satisfies
\begin{equation*}
\Vert A - \widetilde{A} \Vert_2 \le 2 t(r, N, n) t(r, M, m) \sigma_{r+1}(A).
\end{equation*}
\end{theorem}
\begin{proof}
Let $F$ be a difference between $A$ and $Z$. Divide $A$, $Z$
and
$F$ into blocks:
\begin{equation*}
A = \begin{bmatrix}A_{11}&A_{12}\\A_{21}&A_{22}\end{bmatrix}, \;\;
Z = \begin{bmatrix}Z_{11}&Z_{12}\\Z_{21}&Z_{22}\end{bmatrix}, \;\;
F = \begin{bmatrix}F_{11}&F_{12}\\F_{21}&F_{22}\end{bmatrix}.
\end{equation*}
Without the loss of generality, let us assume ``basis'' rows and columns constructed using 
\ref{def:CGR} to be the first
rows and columns of $A$. Introduce singular row-vectors
$V=\begin{bmatrix} W & V_2 \end{bmatrix}$ and a matrix of noise
$E=\begin{bmatrix}E_1 & E_2 \end{bmatrix}$.
So, a CGR-approximation is the following:
\begin{equation*}
\widetilde{A} = \begin{bmatrix}A_{11} \\ A_{21}\end{bmatrix} (USW)^\dagger
    \begin{bmatrix} A_{11} & A_{12} \end{bmatrix}.
\end{equation*}
We have $\begin{bmatrix} A_{11} & A_{12} \end{bmatrix} = USV+E$ and $E$ is
orthogonal to $U$ by the construction:
\begin{equation*}
(USW)^\dagger \begin{bmatrix} A_{11} & A_{12} \end{bmatrix} = \begin{bmatrix}
    W^\dagger W & W^\dagger V_2 \end{bmatrix}.
\end{equation*}
So, we reduce the approximation to:
\begin{equation*}
\widetilde{A} = \begin{bmatrix} A_{11} \\ A_{21} \end{bmatrix} \begin{bmatrix}
W^\dagger W & W^\dagger V_2 \end{bmatrix}.
\end{equation*}
Introduce $C_V$ as the minimum norm solution of $WC_V = V_2$:
\begin{equation*}
C_V = W^\dagger V_2.
\end{equation*}
Then, rewrite $A_{11}$:
\begin{equation*}
A_{11} = USW+E_1,
\end{equation*}
and an approximation of $A_{11}$ reads
\begin{equation*}
\widetilde{A}_{11} = A_{11} W^\dagger W = USW+E_1 W^\dagger W = A_{11}
    -(E_1-E_1 W^\dagger W),
\end{equation*}
and their difference is
\begin{equation*}
\widetilde{A}_{11}-A_{11} = -(E_1-E_1 W^\dagger W).
\end{equation*}
Repeat this process for $A_{12}$:
\begin{equation*}
A_{12} = USV_2+E_2,
\end{equation*}
its approximation:
\begin{equation*}
\widetilde{A}_{12} = A_{11} C_V = USV_2+E_1 C_V = A_{12}+E_1 C_V-E_2,
\end{equation*}
and corresponding difference:
\begin{equation*}
\widetilde{A}_{12}-A_{12} = E_1 C_V - E_2.
\end{equation*}
So, the approximation error of the first $n$ rows is the following:
\begin{equation}
\label{eq:CGR_A11A12}
\begin{bmatrix} \widetilde{A}_{11} & \widetilde{A}_{12} \end{bmatrix}
    - \begin{bmatrix} A_{11} & A_{12} \end{bmatrix} = -E R,
\end{equation}
where
\begin{equation}
\label{eq:CGR_R}
R = \begin{bmatrix} I_{m \times m} - W^\dagger W & -C_V \\ 0_{(M-m) \times m}
    & I_{(M-m) \times (M-m)} \end{bmatrix}.
\end{equation}
To continue with the approximation error of all other rows, we use $C_Z$,
introduced in the previous theorem (Equation \eqref{eq:C_Z}). We have
\begin{equation*}
A_{21} = C_Z(A_{11}-F_{11})+F_{21},
\end{equation*}
its approximation $\widetilde{A}_{21}$:
\begin{multline*}
\widetilde{A}_{21} = A_{21} W^\dagger W = C_Z \widetilde{A}_{11}-(C_Z F_{11}
    -F_{21}) W^\dagger W =\\= C_Z A_{11}-C_Z (E_1 - E_1 W^\dagger W)
    - (C_Z F_{11}-F_{21})W^\dagger W,
\end{multline*}
and its approximation error:
\begin{equation*}
\widetilde{A}_{21}-A_{21} = (C_Z F_{11}-F_{21})(I_{m \times m}-W^\dagger W)
    -C_Z(E_1-E_1 W^\dagger W).
\end{equation*}
And, finally, consider $A_{22}$:
\begin{equation*}
A_{22} = C_Z(A_{12}-F_{12})+F_{22} = C_Z((A_{11}-E_1)C_V+E_2)
    -C_Z F_{12}+F_{22},
\end{equation*}
its approximation $A_{22}$:
\begin{equation*}
\widetilde{A}_{22} = A_{21} C_V = C_Z(A_{11}-F_{11})C_V+F_{21}C_V,
\end{equation*}
and its approximation error:
\begin{equation*}
\widetilde{A}_{22}-A_{22} = -(C_ZF_{11}-F_{21})C_V+(C_Z F_{12}-F_{22})
    +C_Z E_1 C_V-C_Z E_2.
\end{equation*}
So, the error of approximation of $\begin{bmatrix} A_{21} & A_{22}
\end{bmatrix}$
is the following:
\begin{equation}
\label{eq:CGR_A21A22}
\begin{bmatrix} \widetilde{A}_{21} & \widetilde{A}_{22} \end{bmatrix}
    - \begin{bmatrix} A_{21} & A_{22} \end{bmatrix} = \begin{bmatrix} C_Z &
    -I_{(N-n) \times (N-n)}\end{bmatrix} F R - C_Z E R,
\end{equation}
where $R$ was defined earlier in \eqref{eq:CGR_R}. Combining
\eqref{eq:CGR_A11A12} and \eqref{eq:CGR_A21A22} we get the total approximation
error:
\begin{equation}
\label{eq:CGR_error}
\widetilde{A}-A = LFR-PER,
\end{equation}
where $L$ and $P$ are defined as
\begin{equation*}
L = \begin{bmatrix} 0_{n \times n} & 0_{n \times (N-n)} \\ C_Z &
-I_{(N-n) \times (N-n)}\end{bmatrix}, \;
P = \begin{bmatrix} I_{n \times n} \\ C_Z \end{bmatrix}.
\end{equation*}
As $Z$ is the best rank-$r$ approximation of $A$ and $USV$ is the best rank-$r$
approximation of $\begin{bmatrix} A_{11} & A_{12}\end{bmatrix}$, spectral norms
of $F$ and $E$ are upper-bounded by $(r+1)$-th singular value of $A$:
\begin{equation}
\label{eq:EFnorm}
\Vert E \Vert_2 \le \Vert F \Vert_2 = \sigma_{r+1}(A).
\end{equation}
Spectral norms of $L$ and $R$ were already discussed in Theorem
\ref{thm:pseudoapproximation}, but
now $R$ is based on an $m$-by-$(M-m)$ matrix $C_V$, built on the rank-$r$ matrix
$V$, so
\begin{equation*}
\Vert L \Vert_2 = \Vert P \Vert_2 = \sqrt{1+\Vert C_Z \Vert_2^2} \le
t(r, N, n),
\end{equation*}
\begin{equation*}
\Vert R \Vert_2 = \sqrt{1+\Vert C_V \Vert_2^2} \le t(r, M, m),
\end{equation*}
By combining Equation \eqref{eq:CGR_error} with upper bounds on spectral norms
of each matrix, we finish the proof for an estimation in the spectral norm.
\end{proof}

Just like $t(r, N)$, practical values and theoretical estimations of
$t(r, N, n)$ can be very different. 
One of possible ways to solve this discrepancy is to use the Chebyshev norm
instead of the spectral norm. We give here the estimates from
\cite{zo-psa-2016,zo-psa-2017pre} and then propose an alternative way to derive a similar 
estimate based on Theorem \ref{thm:CGRapproximation}.

\begin{theorem}[CGR-approximation error in Chebyshev norm
\cite{zo-psa-2016,zo-psa-2017pre}]
\label{thm:CGRapproximation-zo}
Let $A$ be an $N$-by-$M$ complex or real matrix. Then, there exists a
CGR-approximation $\widetilde{A} = CGR$, based on $n$ rows (forming matrix $R$,
where $r \le n \\N$), $m$ columns (forming matrix $C$, where $r \le m \ll M$)
and a kernel matrix $G$ such that
\begin{equation*}
\Vert A - \widetilde{A} \Vert_C \le \sqrt{\frac{(n+1)(m+1)}{(n+1-r)(m+1-r)}}
\sigma_{r+1}(A).
\end{equation*}
\end{theorem}

The proof of this theorem is based on the properties of submatrices of the maximal
\textit{
projective volume}, which is  the multiplication of leading
singular values (instead of all singular values in the case of our rectangular
volume). If suboptimal submatrices are used, the estimate holds with an additional factor.
In the next Theorem we propose a method how to find a sub-optimal
submatrix to build a CGR-approximation and prove the same bounds as in \cite{zo-psa-2016,zo-psa-2017pre}
with an additional factor not larger than $2$.

\begin{theorem}[Rectangular CGR-approximation error in Chebyshev norm]
\label{thm:CGRapproximation2}
Let $A$ be an $N$-by-$M$ complex or real matrix, matrix $Z$ be its
best
rank-$r$ approximation in spectral and Frobenius norms. Then, a rectangular
CGR-approximant $\widetilde{A}$ (from Definition \ref{def:CGR}), based on
$n$
rows ($r \le n \ll N$) and $m$ columns ($r \le m \ll M$), 
satisfies
\begin{equation*}
\Vert A - \widetilde{A} \Vert_C < 2 \sqrt{\frac{(n+1)(m+1)}
    {(n+1-r)(m+1-r)}} \sigma_{r+1}(A).
\end{equation*}
\end{theorem}
\begin{proof}
Since we use the same approximant, as in Theorem
\ref{thm:CGRapproximation}, we
can reuse equations \eqref{eq:CGR_A11A12} and \eqref{eq:CGR_A21A22}. Obviously,
the error of approximation of the first
$n$ rows of $A$ (Equation \eqref{eq:CGR_A11A12}) in the Chebyshev norm is less,
than a corresponding error of all other rows of $A$ (Equation
\eqref{eq:CGR_A21A22}). So,the  norm of the error is bounded by the
following
inequality:
\begin{equation}
\label{eq:chebestimation}
\Vert \widetilde{A}-A \Vert_C \le \Vert LFR \Vert_C + \Vert C_ZER \Vert_C,
\end{equation}
with matrices $L$, $R$ and $P$ defined in Theorem
\ref{thm:CGRapproximation}. From Remark \ref{remark:basis} about ``basis''
selection for estimations in the Chebyshev norm, we know the following:
\begin{equation*}
\forall i: \: \Vert (C_Z)_{i,:} \Vert_2^2 \le \frac{r}{n+1-r} \Longrightarrow
\forall i: \: \Vert (C_Z)_{i,:} \Vert_2^2 < \frac{n+1}{n+1-r},
\end{equation*}
\begin{equation*}
\forall i: \: \Vert (C_Z)_{i,:} \Vert_2^2 \le \frac{r}{n+1-r} \Longrightarrow
\forall i: \: \Vert L_{i,:} \Vert_2^2 \le \frac{n+1}{n+1-r},
\end{equation*}
\begin{equation*}
\forall j: \: \Vert (C_V)_{:,j} \Vert_2^2 \le \frac{r}{m+1-r} \Longrightarrow
\forall j: \: \Vert R_{:,j} \Vert_2^2 \le \frac{m+1}{m+1-r}.
\end{equation*}
So, we get bounds for both summands of \eqref{eq:chebestimation}:
\begin{equation*}
\Vert LFR \Vert_C \le \sqrt{\frac{(n+1)(m+1)}{(n+1-r)(m+1-r)}} \Vert F \Vert_2,
\end{equation*}
\begin{equation*}
\Vert C_ZER \Vert_C < \sqrt{\frac{(n+1)(m+1)}{(n+1-r)(m+1-r)}} \Vert E \Vert_2.
\end{equation*}
Using bounds on spectral norms of $E$ and $F$ from Equation \eqref{eq:EFnorm}, we
finish the proof.
\end{proof}

\begin{remark}
\label{rmk:norm}
As can be seen from the proof of Theorem~\ref{thm:CGRapproximation2},
the Chebyshev
norm of an approximation error practically linearly depends on the maximum
per-row
Euclidean norm of matrices $C_Z$ and $C_V^*$.
\end{remark}

Latter remark explains one of the possible ways to select ``basis'' rows and
columns to
construct an approximant from the definition \ref{def:CGR} constructively.
We choose these rows in such a way that
the maximum per-row (per-column) Euclidean norm
of a matrix of coefficients should be as small as possible.

\subsection{Rectangular maximal volume algorithm}
\label{sec:rect_maxvol}
As it follows from Theorem
\ref{thm:CGRapproximation2} and Remark \ref{rmk:norm}, one of the practical
ways to
reduce an approximation error in the Chebyshev norm is to select such ``basis''
rows
and
columns, that corresponding minimum norm solutions $C_Z$ and $C_V$ have
small upper bounds on per-row or per-column Euclidean norm. Without the loss of
generality, we reduced initial problem (of building better approximation) to
decreasing the maximum per-row norm of the matrix $C_Z$ from Equation
\eqref{eq:C_Z}.
Let us formalize this smaller problem: we have a $N \times r$ real or complex
matrix
$A$ (with $N \ge r$) and we need to find such a $K \times r$ submatrix
$\widehat{A}$ with a complementary $(N-K) \times r$ submatrix $\widetilde{A}$,
that the minimum norm solution $\widetilde{C}$ of equation
\begin{equation*}
\widetilde{C} \widehat{A} = \widetilde{A},
\end{equation*}
has the minimal possible upper bound on Euclidean length of each row.

We propose an iterative greedy maximization of the volume of $\widehat{A}$ by
an
extension by a single row on each iteration. Let us show that it is equal to
the greedy minimization of the maximum per-row norm of $\widetilde{C}$. Let us
assume we
already have preselected a submatrix $\widehat{A}$ and extend it
with an $i$-th row of $\widetilde{A}$. Then, the rectangular volume of the
extended
$\widehat{A}$ will increase by a factor of $\sqrt{1+\Vert \widetilde{C}_i
\Vert_2^2}$ due to Equation \eqref{eq:volume_grow} from Theorem
\ref{rectmaxvol:est}. A greedy maximization of the volume of $\widehat{A}$
simply means we select the row of the maximum length from $\widetilde{C}$.
\emph{So, a
greedy reduction of an upper bound on a per-row norm of $\widetilde{C}$ is the
same,
as a greedy maximization of the volume of $\widehat{A}$}. An iterative greedy
maximization of the rectangular volume is very similar to the Dykstra algorithm
\cite{dykstra1971augmentation} for an optimal experimental design.

Suppose we already have a good $M \times r$ submatrix $\widehat{A}$ with $K > M
\ge r$ and linearly independent columns and add the $i$-th row $A_i$ of the
matrix $A$:
\begin{equation*}
\widehat{A} \leftarrow \begin{bmatrix} \widehat{A} \\ A_i \end{bmatrix}.
\end{equation*}
Let the matrix of coefficients $C$ be the minimum norm solution of $A = C
\widehat{A}$:
\begin{equation*}
C = A \widehat{A}^\dagger.
\end{equation*}
This means that we have to recompute $C$:
\begin{equation*}
C \leftarrow A \begin{bmatrix}\widehat{A} \\ A_i\end{bmatrix}^\dagger.
\end{equation*}
Let $C_i$ correspond to the $i$-th row of $C$. Then,
\begin{equation*}
A \begin{bmatrix} \widehat{A} \\A_i \end{bmatrix}^\dagger = A \begin{bmatrix}
\widehat{A} \\ C_i \widehat{A} \end{bmatrix}^\dagger = A \widehat{A}^\dagger
\begin{bmatrix}I_{M \times M} \\ C_i \end{bmatrix}^\dagger,
\end{equation*}
and
\begin{equation*}
C \leftarrow C \begin{bmatrix} I_{M \times M} \\C_i \end{bmatrix}^\dagger.
\end{equation*}
The pseudo-inverse of $\begin{bmatrix} I_{M \times M} \\ C_i \end{bmatrix}$ can
be obtained via the following formula:
\begin{multline*}
\begin{bmatrix} I_{M \times M} \\ C_i \end{bmatrix}^\dagger = \left(
\begin{bmatrix} I_{M \times M} \\ C_i \end{bmatrix}^* \begin{bmatrix}
I_{M \times M} \\ C_i \end{bmatrix} \right)^{-1} \begin{bmatrix} I_{M \times M}
\\ C_i \end{bmatrix}^* = \\ =\left( I_{M \times M} +C_i^*C_i\right)^{-1}
\begin{bmatrix} I_{M \times M} \\ C_i \end{bmatrix}^*.
\end{multline*}
The inverse of the matrix $I_{M \times M}+C_i^*C_i$ can be computed in a fast
and simple way with the Sherman-Woodbury-Morrison formula:
\begin{equation*}
\left( I_{M \times M}+C_i^*C_i\right)^{-1} = I_{M \times M}-
\frac{C_i^*C_i}{1+C_iC_i^*}.
\end{equation*}
Finally, we get
\begin{equation*}
C \leftarrow \begin{bmatrix}C-\frac{CC_i^*C_i}{1+C_iC_i^*} & \frac{CC_i^*}{1+
C_iC_i^*}\end{bmatrix}.
\end{equation*}
We can also update the squares of lengths of each row of $C$, denoted by a vector
$L$:
\begin{equation*}
\forall j: \: L_j \leftarrow L_j-\frac{| C_j C_i^* |^2}{1+C_i C_i^*}.
\end{equation*}

As can be seen, augmenting the matrix $\widehat{A}$ by a single row of $A$
requires a rank-1 update of $C$. Since the matrix $C$ has $N$ rows and $M$
columns, the update costs $\approx 4NM$ operations (the computation of $CC_i^*$
and a rank-1 update of $C$). So, each addition of a row to $\widehat{A}$ is
similar to the iteration of the original \textbf{maxvol} algorithm, where all
computations inside one iteration are reduced to rank-1 updates.

So we get a very simple greedy method, which is formalized in Algorithm
~\ref{alg:2maxvol}. We start from a non-singular square submatrix $\widehat{A}$
and corresponding $C$ (we get them from \textbf{maxvol} algorithm), then
iteratively add a row to $\widehat{A}$, corresponding to the row of the
maximal length in $C$, recompute $C$ and update the vector $L$ of lengths of
each
row of the matrix $C$. We call this algorithm the \textbf{rect\_maxvol}
algorithm
as it is a natural extension of the original \textbf{maxvol} algorithm for
rectangular submatrices. Iterations can be stopped when a length of each
row of $C$ is less than a given parameter $\tau$, assuring that the multiplier
for
the approximation error in the Chebyshev norm will not be higher, than
$\sqrt{1+\tau^2}$.
\begin{algorithm}[H]
\caption{\textbf{rect\_maxvol} (``Greedy'' maximization of the volume of
submatrix)}
\label{alg:2maxvol}
\begin{algorithmic}[1]
\REQUIRE A full-rank $A \in \mathbb{C}^{N \times r}, N > r,$ parameter $\tau$
\ENSURE A submatrix $\widehat{A}$, a set of pivot rows $\{p\}$ and a matrix of
coefficients $C$ such, that $A = C \widehat{A}, \forall i \notin \{p\}:
\:\Vert C_i \Vert_2 \le \tau$
\STATE Start with a non-singular square submatrix $\widehat{A}$
\COMMENT{Result of the \textbf{maxvol}}
\STATE $\{p\} \leftarrow$ pivot rows$;\; C \leftarrow A \widehat{A}^{-1};\;
\forall i:\: L_i \leftarrow \Vert C_i \Vert_2^2$ \COMMENT{Result of the
\textbf{maxvol}}
\STATE $i \leftarrow \mathrm{argmax}_{i \notin \{p\}}(L_i)$ \COMMENT{Find
maximal row in $C$}
\WHILE{$L_i > \tau^2$}
\STATE $\{p\} \leftarrow \{p\} + i$ \COMMENT{Extend set of pivots}
\STATE $\widehat{A} \leftarrow \begin{bmatrix} \widehat{A} \\ A_i\end{bmatrix}$
\COMMENT{Extend $\widehat{A}$}
\STATE $C \leftarrow \begin{bmatrix} C - \frac{CC_i^*C_i}{1+C_iC_i^*} &
\frac{CC_i^*}{1+C_iC_i^*} \end{bmatrix}$
\COMMENT{Rank-1 update of $C$}
\STATE $\forall j:\: L_j \leftarrow L_j-\frac{|C_jC_i^*|^2}{1+C_iC_i^*}$
\COMMENT{Update lengths of rows of $C$}
\STATE $i \leftarrow \mathrm{argmax}_{i \notin \{p\}}(L_i)$ \COMMENT{Find
maximal row in $C$}
\ENDWHILE
\IF{$\widehat{C}$ is required to be identity}
\STATE $\widehat{C} = I$
\ENDIF
\RETURN $C, \widehat{A}, \{p\}$
\end{algorithmic}
\end{algorithm}

Numerical experiments with randomly generated $N \times r$ matrices (not
presented here) have shown that Algorithm~\ref{alg:2maxvol} requires only $K
\approx 1.2r$ rows to reach the upper bound of $2.0$ for the length of each
non-``basis'' row of $C$ and only $K \approx 2r$ to reach the upper bound $1.0$
for the length of each non-``basis'' row of $C$. These results are consistent
with the theory from Section \ref{sec:bounds}.

Since we already evaluated the computational cost for the recomputation of $C$,
it is easy to calculate the number of operations, required for
Algorithm~\ref{alg:2maxvol}.
Computation of a non-singular submatrix with a help of the LU decomposition
with
pivoting requires $\mathcal{O}(Nr^2)$ operations. Since operations,
required to compute $CC_i^*$, are already taken into the account in the computation
of $C$, we can get the total complexity of Algorithm~\ref{alg:2maxvol}: it is
$\mathcal{O}(N(2K^2-r^2))$ operations. For the parameter $\tau=1.0$, the
theoretical estimate of $K=2r-1$ gives us the following result: computational
complexity of the Algorithm~\ref{alg:2maxvol} is $\mathcal{O}(Nr^2)$ operations.

\section{Numerical examples}
\label{sec:numerical}
The \textbf{rect\_maxvol} algorithm was implemented in Python (with
acceleration by Cython) and is
available online at \url{https://bitbucket.org/muxas/maxvolpy}. We test the
efficiency of the \textbf{rect\_maxvol} algorithm compared to the
\textbf{maxvol} algorithm on three different applications. 

\subsection{Finding maximum in modulus element in matrix}

\begin{figure}[H]
\centering
\resizebox{12cm}{!}{\input{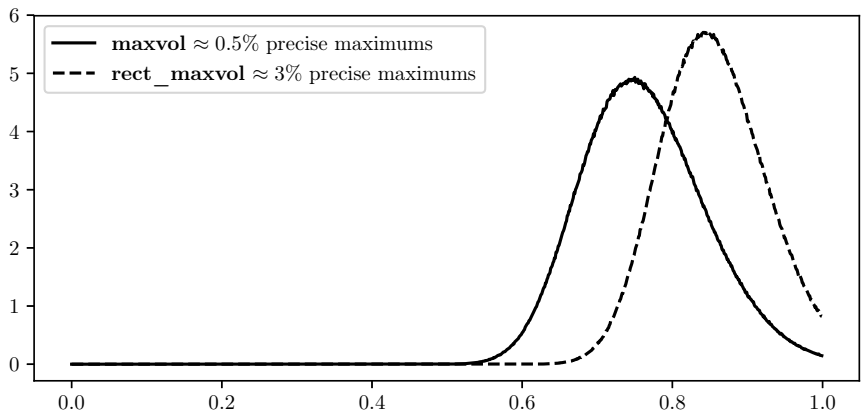}}
\caption{Distribution of the ratio of results of \textbf{maxvol} and
\textbf{rect\_maxvol} over the true maximums for 8163272 random experiments.}
\label{pic:ratiodistribution}
\end{figure}

The \textbf{maxvol} algorithm is a heuristic procedure to find the maximal in
modulus element in a low-rank matrix. We repeat the corresponding experiment
from \cite{gostz-maxvol-2010}. We generate random low-rank matrices as a
multiplication of $3$ matrices,
\begin{equation*}
A = UDV^T,
\end{equation*}
where $U \in \mathbb{R}^{10000 \times 10}$ and $V \in \mathbb{R}^{10000 \times
10}$ are Q-factors of the QR factorization of randomly generated matrices with
uniformly distributed in the interval $[0;1]$ elements and $D \in \mathbb{R}^{
10 \times 10}$ is a randomly generated diagonal matrix with uniformly
distributed in $[0;1]$ elements. Assuming we have a low-rank approximation of
each test matrix, we find the maximal-volume rows and columns of $U$ and $V^T$
correspondingly and measure a ratio of the maximal absolute element on the
intersection of found rows and columns to the maximal absolute element in the
entire matrix. We have measured the latter ratio for each test matrix with two
different ways of finding the maximal-volume rows/columns: by \textbf{maxvol}
and
by \textbf{rect\_maxvol}. Results are presented in
\figurename~\ref{pic:ratiodistribution}.

\subsection{Preconditioning of overdetermined systems}

This example was inspired by the paper \cite{arioli2015preconditioning}, 
where among other techniques the authors have used row selection based on
\textbf{maxvol} algorithm for the preconditioning of a least squares problem.
Here
we show, that the condition number can be made much better using the
\textbf{rect\_maxvol} algorithm.

Assume we need to solve an overdetermined system $Ax=b, A \in \mathbb{C}^{N
\times
r}, r \ll N$, in the least-squares sense:
\begin{equation}
x = \mathrm{argmin} \Vert Ax-b \Vert_2^2.
\end{equation}
This is equivalent to the solution of the normal equations
\begin{equation*}
A^* A x = A^*b.
\end{equation*}

With a help of the \textbf{rect\_maxvol} algorithm we can find the submatrix
$\widehat{A}$ and the matrix of coefficients $C$ such that
\begin{equation*}
A = C \widehat{A}, \;\; C = A \widehat{A}^\dagger, \;\; \widehat{A} \in
\mathbb{C}^{K \times r}, \; C \in \mathbb{C}^{N \times K}, \; K \ge r.
\end{equation*}
Thus, there is a permutation matrix $P \in \mathbb{R}^{N \times N}$ such that:
\begin{equation*}
PA = \begin{bmatrix} \widehat{A} \\ B \end{bmatrix}.
\end{equation*}
Using this partitioning of $A$ into a basic part $\widehat{A}$ and a non-basic
part $B$, we rewrite the residual vector $r$ and the right hand side $b$ as
\begin{equation*}
Pr = \begin{bmatrix} r_{ \widehat{A}} \\ r_B \end{bmatrix},\; Pb =
\begin{bmatrix} b_{ \widehat{A}} \\ b_B \end{bmatrix}.
\end{equation*}

If $x$ is a solution of the system $A^* A x = A^* b$, then $x$ is a solution of
\begin{equation}
\label{eq:solution2}
C^* C \widehat{A}x = C^*b.
\end{equation}
Then, we construct an augmented system
\begin{equation}
\label{eq:finalsystem}
Z \begin{bmatrix}r_B \\ \widehat{A} x \end{bmatrix} = \begin{bmatrix} b_B \\ -
\widehat{C} b_{\widehat{A}} \end{bmatrix}, \;\;\; Z = \begin{bmatrix} I_{
(N-K) \times (N-K)}
& \widetilde{C} \\ \widetilde{C}^* & -I_{K \times K} \end{bmatrix},
\end{equation}
where $\widehat{C}$ is a basic part of $C$ (such that $\widehat{A} =
\widehat{C} \widehat{A}$) and $\widetilde{C} = B \widehat{A}^\dagger$.
If we eliminate the first $(N-K)$ variables, we will get the equation
\eqref{eq:solution2}.
A solution of the system \eqref{eq:finalsystem} consists of 2 parts: solve a
system with the matrix $Z$ and solve a least squares problem with the matrix
$\widehat{A}$. However, if we solve the system with the matrix $Z$ precisely
and put $\widehat{C} = \widehat{A} \widehat{A}^\dagger$, least squares $K
\times r$ problem with the matrix $\widehat{A}$ has unique solution and
can be reduced to $r \times r$ system by finding a good square submatrix in
$\widehat{A}$.

In \cite{arioli2015preconditioning} it was shown that the condition number of
the system \eqref{eq:finalsystem} is the following:
\begin{equation*}
\mathrm{cond}(Z) = \sqrt{1+\Vert \widetilde{C} \Vert_2^2}.
\end{equation*}
Therefore, the condition number of $Z$ is equal to the spectral norm of the
matrix $C$ due to \eqref{eq:topsingval}:
\begin{equation*}
\mathrm{cond}(Z) = \Vert C \Vert_2 \le t(r, N, K)
\end{equation*}
and is bounded only by $t(r, N, K)$.

For experiments, we used $3$ ill-conditioned sparse matrices, available on the
Web: \textbf{illc1850}, \textbf{lp\_osa\_07} and \textbf{Kemelmacher}. In these
model experiments we did not use the sparsity of those matrices, since our goal
was to estimate the final condition number. Efficient implementation of the
\textbf{rect\_maxvol} algorithm for sparse matrices is a topic of ongoing work.
In the Table~\ref{table:preconditioning} we present results of experiments.

\begin{table}[H]
\caption{Comparison of preconditioning by \textbf{maxvol} and
\textbf{rect\_maxvol} algorithms. Time is measured in seconds, rows
corresponds to parameter $K$.}
\label{table:preconditioning}
\begin{tabular}{c|ccc|ccc}
Matrix & \multicolumn{3}{c|}{maxvol} & \multicolumn{3}{c}{rect\_maxvol} \\
name & time & rows & $\Vert C \Vert_2$ & time & rows &
$\Vert C \Vert_2$ \\ \hline
illc1850 & 0.39 & 712 & 15.96 & 0.51 & 1095 & 4.37 \\
lp\_osa\_07 & 3.22 & 1118 & 184.8 & 92.7 & 2184 & 11.66 \\
Kemelmacher & 339.82 & 9693 & 60.93 & 4135.34 & 15237 & 5.17 \\
\end{tabular}
\end{table}

\subsection{Recommender systems}

Another application comes from the field of recommender systems. A
collaborative filtering deals with the \emph{user-product matrix} $A$, which
encodes the ratings for a particular user. The SVD is often used for the
prediction of ratings the user will give to a particular problem. The
\emph{cold start problem} is the problem of the rating for a new user. One of
possible solutions relies on the extremal submatrices. In
\cite{liu2011wisdom} authors proposed to use \textbf{maxvol} to find
\emph{representative users}. This type of factorization is based on a
skeleton approximation of the matrix $A$ using its rows $R$ and its columns
$C$:
\begin{equation*}
A \approx  C \widehat{A}^{-1} R, \;\; A \in \mathbb{C}^{N \times M}, \; C \in
\mathbb{C}^{N \times r},\; R \in \mathbb{C}^{r \times M}, \; \widehat{A} \in
\mathbb{C}^{r \times r},
\end{equation*}
\begin{equation*}
\mathrm{rank}(A) = r \ll \min(N, M),
\end{equation*}
where $\widehat{A}$ is a submatrix of $A$ on the intersection of rows $R$
and columns $C$. At a preprocessing step, the user-product matrix is
approximated by its best low-rank approximation computed by the SVD. Once
columns $C$ or rows $R$ are selected, we can compute weights $X_C$ or $X_R$
from the least squares approximation:
\begin{equation*}
A \approx C X_C \approx X_R R.
\end{equation*}
This decomposition has a very simple meaning: ratings of all products for any
given user is a linear combination of ratings of the ``most representative
users''
and ratings, given by all users, of any given product is a linear combination
of
ratings of the ``most representative products''. When new user appears in such
a
database, he/she can be asked to rank the ``most representative products'' to
update the decomposition. On the other hand, when the new product is added, the
``most representative users'' can be asked to rank it to update the
decomposition.

We applied the \textbf{rect\_maxvol} algorithm to choose representative
users or items and construct the corresponding approximation. For numerical
examples we used the MovieLens dataset
\url{http://grouplens.org/datasets/movielens/} with $10$ million ratings with
$10000$ movies by $72000$ users. At first, we computed the best rank-$k$
approximation from the SVD. Then, we computed either \textbf{maxvol} or
\textbf{rect\_maxvol} rows/columns. To measure the quality, we used the
\emph{coverage}, \emph{diversity} and \emph{precision} criterias, same as in
\cite{liu2011wisdom}:

\begin{itemize}
\item \textbf{Coverage}: proportion of users (movies) which rated
 (were rated by) any of the representative movies (users),
\item \textbf{Diversity}: proportion of users (movies) which rated
 (were rated by) any, but less than 10 \% of the representative movies (users),
\item \textbf{Precision}: proportion of good recommendations among
 the top k recommendations.
\end{itemize}
Each metric was calculated as an average for every user (movie). Corresponding
results are shown in \tablename~\ref{table:coverage} and
\tablename~\ref{table:precision}.

\begin{table}[H]
\caption{Coverage and diversity of \textbf{maxvol} and \textbf{rect\_maxvol}
representatives.}
\label{table:coverage}
\begin{tabular}{c|ccc|ccc}
 & \multicolumn{3}{c|}{user} & \multicolumn{3}{c}{movie} \\
 & k & coverage & diversity & k & coverage & diversity \\ \hline
\textbf{maxvol} & 100 & 0.89 & 0.6 & 20 & 0.94 & 0.11 \\
\textbf{rect\_maxvol} & 50 & 0.89 & 0.6 & 10 & 0.91 & 0.14
\end{tabular}
\end{table}

\begin{table}[H]
\caption{Precision at 10 for 5 derivatives from MovieLens data}
\label{table:precision}
\begin{tabular}{c|c|c|ccccc}
\multirow{2}{*}{Type} & \multirow{2}{*}{k} & \multirow{2}{*}{criteria} &
\multicolumn{5}{|c}{Dataset} \\
 & & & 1 & 2 & 3 & 4 & 5 \\ \hline
\multirow{2}{*}{\textbf{maxvol}} & \multirow{2}{*}{20} & Precision at
10 & 0.46 & 0.45 & 0.47 & 0.45 & 0.46 \\
 & & representative movies & 20 & 20 & 20 & 20 & 20 \\ \hline
\multirow{2}{*}{\textbf{rect\_maxvol}} & \multirow{2}{*}{10} &
Precision at 10 & 0.5 & 0.49 & 0.52 & 0.49 & 0.5 \\
 & & representative movies & 15 & 14 & 16 & 14 & 15
\end{tabular}
\end{table}

It is very interesting, that it is better to select $20$ rows using the best
rank-$10$ approximation, rather than compute the best rank-$20$ approximation
with the classical \textbf{maxvol} algorithm. This should definitely be studied
in more
details. 

\section{Related work}
\label{sec:optdesign}

Related theoretical work is mostly based on estimations of a skeleton-type
approximation error. We cited different papers, where such an estimation is
based on $(r+1)$-th singular value using $r$ rows and $r$
columns for approximation itself. However, multiplier of that $(r+1)$-th
singular value is rather high and, thus, can be reduced by using more, than $r$
rows and columns. Recent papers \cite{zo-psa-2016} and \cite{zo-psa-2017pre} on
this theme show dependency of using additional rows and columns on
an investigated error multiplier.

Algorithmical approaches, similar to described in this paper, were also
provided in an optimal experimental design. The problem of an optimal
experimental
design is based on the following linear regression model:
\begin{equation}
y = A x + n,
\end{equation}
where $y$ is a vector of $N$ responses, $A$ is a $N$-by-$r$ matrix of
independent variables, $x$ is a vector of regression coefficients and $n$ is
a vector of errors. Each variable of $n$ is assumed to be independent and
normally distributed with the same variance. Problem here is to select such a
subset of experiments (rows of $A$ with corresponding $y$), that influence of a
white noise $n$ is as small, as possible. One of possible solutions is to
select such a
submatrix $\widehat{A}$, based on rows of $A$, which minimizes generalized
variance \cite{kiefer1961optimum}, which is equivalent to maximizing $\det
\widehat{A}^* \widehat{A}$. Such an optimization criteria is
usually called \emph{D-optimality} in an optimal experimental design
literature.
So, maximization of the rectangular volume is the same, as \emph{D-optimality}
criteria. Main difference of well-known algorithms of finding good
\emph{D-optimal} submatrices
\cite{mitchell1974algorithm,dykstra1971augmentation,fedorov1972theory} and
\textbf{rect\_maxvol}, proposed in this paper, is that our algorithm is based
on a Gauss elimination to find a good submatrix to start with, whereas
algorithms
from \cite{mitchell1974algorithm,dykstra1971augmentation,fedorov1972theory} use
a random submatrix for this purpose.

\section{Conclusion and future work}

Rectangular submatrices have high potential in different column/row sampling
methods. The rectangular volume maximization leads to an efficient
computational algorithm, proved to be useful not only for approximations of
matrices. A construction, proposed in definition~\ref{def:CGR} may
lead
to a new cross approximation technique, which is a subject for future research.

\section*{Acknowledgements}

Work on the problem setting and numerical examples was supported by Russian
Foundation for Basic Research grant 16-31-60095 mol\_a\_dk. Work on theoretical
estimations of approximation error and the practical algorithm was supported by
Russian Foundation for Basic Research grant 16-31-00351 mol\_a.

\bibliography{bibtex/our,bibtex/algebra,rect_maxvol}

\end{document}